\documentclass[submission%
]{dmtcs}

\usepackage[latin1]{inputenc}
\usepackage{subfigure}
\usepackage{amsmath}
\usepackage{enumerate}
\usepackage{amssymb}
\usepackage[dvipsnames]{xcolor}
\usepackage[normalem]{ulem}

\newcommand{\ol}[1]{\ensuremath{\overline{#1}}} % Overline

\newtheorem{thm}{Theorem}

\newtheorem{lemma}[thm]{Lemma}
\newtheorem{sublemma}[thm]{Sublemma}

\newtheorem{cor}[thm]{Corollary}
\newtheorem{prop}[thm]{Proposition}
\newtheorem{dfn}[thm]{Definition}
\newtheorem{ex*}[thm]{Example}
\newenvironment{ex}{\begin{ex*}\rm}{\end{ex*}}

\newtheorem{citedfn}[thm]{Definition}

\author{Amanda Cameron\thanks{Email: \email{amanda.cameron@qmul.ac.uk}}%
  \and Alex Fink\thanks{Email: \email{a.fink@qmul.ac.uk}}}
\title{A lattice point counting generalisation of the Tutte polynomial}
\address{Queen Mary University of London, UK}
\received{2015-11-16}
\revised{\today}
\accepted{tomorrow}
\begin{document}
\maketitle
\begin{abstract}
\paragraph{Abstract.}
The Tutte polynomial for matroids is not directly applicable to polymatroids. For instance, deletion-contraction properties do not hold. 
We construct a polynomial for polymatroids which behaves similarly to the Tutte polynomial of a matroid, and in fact contains the same information as the Tutte polynomial when we restrict to matroids. This polynomial is constructed using lattice point counts in the Minkowski sum of the base polytope of a polymatroid and scaled copies of the standard simplex. 
We also show that, in the matroid case, 
our polynomial has coefficients of alternating sign,
with a combinatorial interpretation closely tied to the Dawson partition.

\paragraph{R\'esum\'e.} 
Le polyn\^ome de Tutte pour les matro\"ides n'est pas directement applicable aux polymatro\"ides. Par exemple, les propri\'et\'es de suppression et de contraction ne sont pas v\'erifi\'ees. Nous construisons un polyn\^ome pour les polymatro\"ides  qui se comporte de la m\^eme fa\c{c}on que le polyn\^ome de Tutte d'un matro\"ide et qui contient en fait la m\^eme information que le polyn\^ome de Tutte lorsqu'il est restreint aux matro\"ides. Ce polyn\^ome est construit en comptant les points entiers dans la somme de Minkowski du polytope des bases d'un polymatro\"ide et de copies homoth\'etiques du simplexe standard. Nous montrons aussi que dans le cas d'un matro\"ide, notre polyn\^ome a des coefficients \`a signes alternants, avec une interpr\'etation \'etroitement li\'ee \`a la partition de Dawson.
\end{abstract}

\section{Introduction}
\label{sec:in}

The Tutte polynomial, originally formulated for graphs, has been generalised to apply to matroids.  We recommend \cite{Oxley} as a reference for basic terminology of matroids. 
Let $\mathcal P(E)$ be the power set of~$E$.

\begin{dfn}\label{dfn:Tutte}
Let $M=(E,r)$ be a matroid with ground set $E$ and rank function $r:\mathcal{P}(E)\rightarrow \mathbb{Z^+}\cup\{0\}$. The \emph{Tutte polynomial} of $M$ is 
\begin{displaymath}
    T_M(x,y) = \sum_{S\subseteq E} (x-1)^{r(M)-r(S)}(y-1)^{|S|-r(S)}.
\end{displaymath}    
    \end{dfn}
    
This polynomial has a diverse range of applications, from classifying Tutte invariants --- properties of matroids or graphs which can be enumerated by an evaluation of Tutte --- to practical applications in coding theory. A natural extension of matroids are polymatroids, which are a class of objects formed by relaxing one matroid rank axiom.  The Tutte polynomial does not directly apply to these and provide similarly useful results. 

\begin{dfn}
A \emph{polymatroid} can be described by the ground set $E$ and a rank function $r:\mathcal{P}(E)\rightarrow \mathbb{Z^+}\cup\{0\}$ such that, for $X,Y\in\mathcal{P}(E)$, the following conditions hold:
\begin{itemize}
\item[{\rm P1}.] $r(\emptyset)=0$
\item[{\rm P2}.] If $Y\subseteq X$, $r(Y)\leq r(X)$
\item[{\rm P3}.] $r(X\cup Y)+r(X\cap Y)\leq r(X)+r(Y)$
\end{itemize}
\end{dfn}

We will be principally viewing polymatroids as polytopes. Let $E$ be a finite set, which will serve as the ground set of our 
(poly)matroids.
We work in the vector space $\mathbb R^E$.
For a set $U\subseteq E$, let ${\bf e}_U\in\mathbb R^E$ be the indicator vector of $U$,
and abbreviate ${\bf e}_{\{i\}}$ by ${\bf e}_{i}$.
Let $r:\mathcal{P}(E)\rightarrow\mathbb{N}$ be a rank function,
and $M=(E,r)$ the associated polymatroid. Define the following polytope to be the \emph{extended polymatroid} of $r$:
\begin{displaymath}EP(M)=\{{\bf x}\in\mathbb{R}^E \ | \  {\bf x}\geq 0\mbox{ and }
{\bf x}\cdot{\bf e}_U\leq r(U) \ \mathrm{for \ all} \ U\subset E\}.\end{displaymath}
We also have the polymatroid \emph{base polytope} of $r$,
a face of the extended polymatroid:
\[P(M) = EP(M)\cap\{{\bf x}\in\mathbb{R}^E \ | \  {\bf x}\cdot{\bf e}_E= r(E)\}
= \operatorname{conv} \mathcal B_M.\]
Either one of these polytopes contains all the information in the rank function.
In \cite{postnikov} polymatroid base polytopes are dubbed
``generalised permutohedra''.

The presentation of the Tutte polynomial in Definition~\ref{dfn:Tutte}
is given in terms of the \emph{corank-nullity polynomial}:
up to a change of variables, it is the generating function for subsets $S$ of the
ground set by their corank $r(M)-r(S)$ and nullity $|S|-r(S)$.
The corank-nullity polynomial can be defined for polymatroids, but the resulting
function gives no easily accessible information about basis activities,
and is not even a Laurent polynomial in the variables $x$ and~$y$.

This reflects the difference between matroids and polymatroids 
that matroids have a well-behaved theory of \emph{minors} analogous to graph minors:
for each ground set element one can define a deletion and contraction,
and knowing these two determines the matroid.  
The deletion-contraction recurrence for the Tutte polynomial reflects this structure.
While polymatroids fulfil general notions of deletion and contraction, these do not allow for a similar deletion-contraction recurrence if we attempt to construct 
a modified Tutte polynomial in the directly analogous fashion.
It is however possible to salvage some features of this
recurrence in restricted cases: this is done by \cite{whittle}
for polymatroids where singletons have rank at most~2, where 
the corank-nullity polynomial is still universal for a form of
deletion-contraction recurrence.

Another formula for the Tutte polynomial of a graph or a matroid is defined
in terms of \emph{activities} of \emph{bases}. The following three 
definitions of activity for the polymatroid generalisation are from \cite{kalman}.

\begin{dfn}
A vector ${\bf x}\in\mathbb Z^E$ is called a \emph{base} if ${\bf x}\cdot{\bf e}_E=r (E)$ 
and ${\bf x}\cdot{\bf e}_S\leq r (S)$ for all subsets $S\subseteq E$.
\end{dfn}

Let $\mathcal B_M$ be the set of all bases of~$r$.

\begin{dfn}
A \emph{transfer} is possible from $u_1\in E$ to $u_2\in E$ 
in the base ${\bf x}\in \mathcal B_M\cap\mathbb{Z}^E$
if by decreasing the $u_1$-component of ${\bf x}$ by $1$ and increasing its $u_2$-component by $1$ we get another base.
\end{dfn}

\begin{dfn}
Order the elements of $E$ arbitrarily.
\begin{enumerate}[i.]
\item We say that $u\in E$ is \emph{internally active} with respect to the base $x$ if no transfer is possible in~$x$ \emph{from} $u$ to a smaller element of $E$.
\item We say that $u\in E$ is \emph{externally active} with respect to $x$ if no transfer is possible in~$x$ \emph{to} $u$ from a smaller element of $E$.
\end{enumerate}
\end{dfn}

For $x\in \mathcal B_M\cap\mathbb{Z}^E$\!, let the set of internally active elements with respect to $x$ be denoted with $\mathrm{Int}(x)$, and let $\iota(x)=|\mathrm{Int}(x)|$;
likewise, let the set of externally active elements be denoted with $\mathrm{Ext}(x)$
and $\varepsilon(x)=|\mathrm{Ext}(x)|$.
Let $\ol{\iota}(x),\ol{\varepsilon}(x)$ denote 
the respective numbers of inactive elements. When $M$ is a matroid, these numbers provide an alternative formulation of the Tutte polynomial, \begin{displaymath}\label{eq:act}
T(M;x,y)=\sum_{B\in\mathcal{B}_M} x^{\iota(B)}y^{\varepsilon(B)}.
\end{displaymath}
In this case, the following simplified versions of the definition of activity are more commonly used:

\begin{dfn}
Take a matroid $M=(E,r)$, and give $E$ some ordering. Let $B$ be a basis of $M$. 
\begin{enumerate}[i.]
\item We say that $e\in E-B$ is \emph{externally active} with respect to $B$ if $e$ is the smallest element in the unique circuit contained in $B\cup e$, with respect to the ordering on $E$.
\item We say that $e\in B$ is \emph{internally active} with respect to $B$ if $e$ is the smallest element in the unique cocircuit in $(E\setminus B)\cup e$.
\end{enumerate}
\end{dfn}

These definitions are the analogies of those originally formulated using spanning trees of graphs. 

The analogy between the internal and external polynomials of a polymatroid and the same polynomials under the graph definitions suggests that a two-variable polynomial similar to Tutte can indeed be found for polymatroids. In the case of a polymatroid, the definitions are as follows. Note that these polynomials do not depend on the order on~$E$ that was used to define them. 

\begin{citedfn}[{\cite{kalman}}]
\label{activity}
Define the \emph{internal polynomial} and \emph{external polynomial} of $r$ by
\begin{displaymath}I_r(\xi)=\sum_{x\in \mathcal B_M\cap\mathbb{Z}^E} \xi^{\ol{{\iota}}(x)} \quad \mathrm{and} \quad X_r(\eta)=\sum_{x\in \mathcal B_M\cap\mathbb{Z}^E} \eta^{\ol{\varepsilon}(x)}.\end{displaymath}
\end{citedfn}

Our work generalises formula \eqref{eq:act}, creating a two-variable polynomial which is equivalent to Tutte for matroids (Theorem~\ref{lpg})
 and specialises to the two activity polynomials above for polymatroids (Theorem~\ref{thm:activity}).  That is, the invariant we construct is the bivariate analogue of K\'{a}lm\'{a}n's activity polynomials, which is something his paper sought.
Note that the number of bases in $M$ is the number of lattice points in $P(M)$. We will thus form a polynomial which counts the lattice points in a particular polytope which we construct
from $P(M)$ in a way which introduces the stylistically necessary two variables. This will be described in full detail in the next section.

K\'alm\'an's original interest in these objects related to enumerating
spanning trees of bipartite graphs according to their vector of degrees at
the vertices on one side.  In this context \cite{oh} has investigated a polyhedral
construction similar to ours, as a way of proving Stanley's pure O-sequence
conjecture for cotransversal matroids.  \cite{kalman2}
have since explained a duality preserving $I_r(\xi)$ arising from planar maps,
using root polytopes.

This paper is organised as follows: In Section \ref{sec:hints} we begin by explaining the construction of our polymatroid polynomial, followed by how this is directly related to the Tutte polynomial when we restrict the domain to matroids. The main theorem, Theorem \ref{lpg}, gives the Tutte polynomial as a sum of lattice points in a particular polytope. The section ends with properties of our polynomial which hold under polymatroid generalisation. Finally, in Section \ref{sec:coeffs}, we again restrict our attention to matroids. We give a geometric interpretation of the coefficients of our polynomial, by way of a particular subdivision of the relevant polytope.

\section{The polynomial}
\label{sec:hints}

\subsection{Construction}
\label{sec:const}

Note that the number of lattice points in a base polytope is equal to the number of bases in its underlying matroid. We will form a polynomial which counts the lattice points of a particular Minkowski sum of polyhedra.

Let $\Delta$ be the standard simplex in $\mathbb R^E$ of dimension equal to $r(M)-1$,
that is 
\[\Delta = \operatorname{conv}\{\mathbf e_{i}:i\in E\},\]
and $\nabla$ be its reflection through the origin,
$\nabla = \{-x : x\in\Delta\}$.
The faces of $\Delta$ are the polyhedra
\[\Delta_S = \operatorname{conv}\{\mathbf e_{i}:i\in S\}\]
for all nonempty subsets $S$ of~$E$; similarly, 
the faces of $\nabla$ are the polyhedra $\nabla_S$ given as the reflections of the $\Delta_S$.

We consider $P(M)+u\Delta+t\nabla $ where $M=(E,r)$ is any polymatroid and $u,t\in\mathbb{N}$. 
We are interested in the lattice points in this sum, which are the vectors
that can be turned into bases of~$M$ by incrementing a coordinate $t$ times
and decrementing one $u$ times.
If we write the polytope in the form \begin{displaymath}\{{\bf x}\in\mathbb{R}^d \ | \ A{\bf x}\leq {\bf b}, \ A\in\mathbb{Z}^{m\times d},{\bf b}\in\mathbb{Z}^m\}\end{displaymath} where $d$  is the dimension of the polytope and $A{\bf x}\leq {\bf b}$ is a system of rank inequalities describing the polytope, then finding the number of lattice points inside the polytope is equivalent to finding the number of integer solutions ${ \bf x}$ to $A{\bf x}\leq {\bf b}$. 
By Theorem 7 of \cite{mcmullan}, this number of lattice points
is a polynomial in $t$ and~$u$, of degree $\dim(P(M)+u\Delta+t\nabla) = |E|-1$.
This polynomial we rewrite in the form \begin{displaymath}Q(M;t,u):=\#(P(M)+u\Delta+t\nabla )\cap\mathbb{Z}^E=\sum_{i,j} c_{ij}\binom{u}{j}\binom{t}{i}.\end{displaymath}
Changing variables gives the polynomial \begin{displaymath}Q'(M;x,y)=\sum_{ij}c_{ij}(x-1)^i(y-1)^j\end{displaymath}
where the $c_{ij}$ are equal to those in the previous equation.
This change of variables is chosen so that
applying it to $\#(u\Delta_X+t\nabla_Y)$ yields $x^iy^j$,
where $\Delta_X$ and $\nabla_Y$ are faces of $\Delta$ and $\nabla$ 
of respective dimensions $i$ and~$j$.  This will allow for a
combinatorial interpretation of its coefficients in Theorem~\ref{coeffs}.

One motivation for this particular Minkowski sum is that it provides a polyhedral
translation of K\'alm\'an's construction of activities in a polymatroid.

\begin{lemma}
\label{old}Let $P$ be a polymatroid polytope. At every point $f\in P$, attach the scaled simplex \[f+t\operatorname{conv}(\{-e_i \ | \ i \ \textrm{is \ internally \ active \ in} \ f \ \textrm{or} \ i\notin f \}).\] This operation partitions $P+t\nabla $
into a collection of translates of faces of $t\nabla $,
with the simplex attached at $f$ having codimension $\ol{\iota}(f)$
within $P$.
\end{lemma}

The following is a direct consequence of this lemma,
and its exterior analogue which arises from replacing $\iota$ and $\nabla$
with $\varepsilon$ and $\Delta$.  

\begin{thm}\label{thm:activity}
Let $M$ be a polymatroid with rank function $r$ and ground set $E$. Then
\[I_r(\xi) = \xi^{|E|-1}Q'(M;\frac{1}{\xi},1) \ \textrm{and} \ X_r(\eta) = \eta^{|E|-1}Q'(M;1,\frac{1}{\eta}).\]
\end{thm}

It is this result which first motivated
the particular change of basis made from $Q$ to~$Q'$,
since an $i$-dimensional face of $t\Delta$ 
has $\displaystyle\binom{t+i}i = \sum_{k=0}^i\binom ik\binom ti$ lattice points.

\subsection{Relation to Tutte}
\label{tutte}

When we restrict $M$ to be a matroid but not a polymatroid, 
$Q'(M;x,y)$ is an evaluation of the Tutte polynomial, and in fact one that contains
precisely the same information.  
As such, the Tutte polynomial can be evaluated by lattice point counting methods.  The main theorem of this section is Theorem \ref{lpg}, but we will first present the relationship between $Q'(M;x,y)$ and Tutte.

\begin{thm}\label{thm:T to Q}
Let $M$ be a matroid. Then we have that
\begin{equation}
Q'(M;x,y)=\dfrac{x^{|E|-r}y^r}{x+y-1}T(\dfrac{x+y-1}{y},\dfrac{x+y-1}{x})
\end{equation}
\end{thm}

We can invert this formula by setting $x'=\dfrac{x+y-1}{y},y'=\dfrac{x+y-1}{x}$, rearranging, and then relabelling. 
\begin{thm}\label{thm:Q to T}
Let $M$ be a matroid. Then
\begin{equation}
T(M;x,y)=-\,\dfrac{(xy-x-y)^{|E|-1}}{(-y)^{r-1}(-x)^{|E|-r-1}}Q'(\dfrac{-x}{xy-x-y},\dfrac{-y}{xy-x-y})
\end{equation}
\end{thm}

We conjecture that there is a relationship between
our formula for the Tutte polynomial and
the algebro-geometric formula for the Tutte polynomial in \cite{speyer},
given that the computations on the Grassmannian in that work are done 
in terms of~$P(M)$, the moment polytope of a certain torus orbit closure,
and that $\Delta$ and $\nabla$ are the moment polytopes of the 
two dual copies of~$\mathbb P^{n-1}$, the $K$-theory ring of whose product
$\mathbb Z[x,y]/(x^n,y^n)$ is identified with the ambient ring of the Tutte polynomial.

\begin{ex}\label{ex:1}
Let $M$ be the matroid on ground set $[3]=\{1,2,3\}$ with
$\mathcal B_M = \{\{1\},\{2\}\}$.
When $t=2$ and $u=1$, the sum $P(M)+u\Delta+t\nabla $ is the polytope of Figure~\ref{fig:1}, with 16 lattice points.
\begin{figure}
\begin{center}
\includegraphics[scale=1]{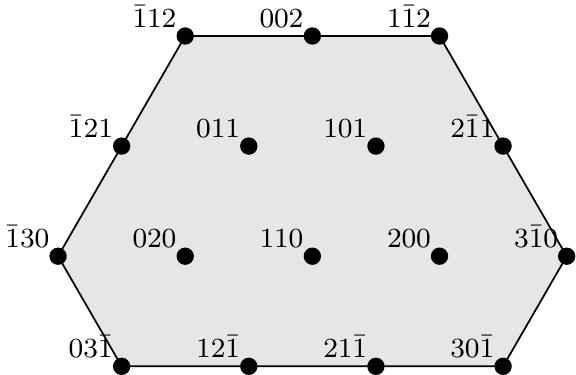}
\end{center}
\caption{The polytope $P(M)+u\Delta+t\nabla $ of Example~\ref{ex:1}.
The coordinates are written without parentheses or commas,
and $\bar1$ means $-1$.}\label{fig:1}
\end{figure}

To compute $Q(M;x,y)$, it is enough to count the lattice points in
$P(M)+u\Delta+t\nabla $ for a range of $u$ and $t$. 
Since $Q$ is a polynomial of degree~2, it is sufficient to take
$t$ and $u$ nonnegative integers with sum at most~2.
These are the black entries in the table below:
\begin{center}
\begin{tabular}{l|ccc}
$t$ $\backslash$ $u$ & 0 & 1 & 2 \\\hline
0 & 2 & 5 & 9 \\
1 & 5 & 10 & \textcolor{Gray}{16} \\
2 & 9 & \textcolor{Gray}{16} & \textcolor{Gray}{24}
\end{tabular}
\end{center}
We can then fit a polynomial to this data, and find
\[Q(M;t,u) = \binom t2 + 2tu + \binom u2 + 3t + 3u + 2,\]
so
\begin{align*}
Q'(M;x,y) &= (x-1)^2 + 2(x-1)(y-1) + (y-1)^2 + 3(x-1) + 3(y-1) + 2
\\&= x^2 + 2xy + y^2 - x - y.
\end{align*}
Finally, by Theorem~\ref{thm:Q to T},
\begin{align*}
T(M;x,y) &= -\,\frac{(xy-x-y)^2}{(-y)^0(-x)^1}
\left(\frac{y^2 + 2xy + x^2}{(xy-x-y)^2}+\frac{y + x}{xy-x-y}\right)
\\&= xy + y^2
\end{align*}
which is indeed the Tutte polynomial of~$M$.
\end{ex}

We can use Theorem \ref{thm:Q to T} to give a formula for the Tutte polynomial directly in terms of lattice point counting, as follows:

\begin{thm}
\label{lpg}
Let $M=(E,r)$ be a matroid with $r(M)=r$. Then
\[T(x,y)=(xy-x-y)^{|E|}(-x)^r(-y)^{|E|-r}\sum_{u,t\geq 0} Q(t,u)\cdot(\dfrac{y-xy}{xy-x-y})^t\cdot(\dfrac{x-xy}{xy-x-y})^u\]
\end{thm}

A further substitution and simple rearrangement gives the following corollary, included for the sake of completeness.

\begin{cor}
Let $M=(E,r)$ be a matroid with $r(M)=r$. Then
\[\sum_{u,t\geq 0}Q(t,u)v^tw^u=\dfrac{1}{(1-v)^{|E|-r}(1-w)^{r}(1-vw)}\cdot T(\dfrac{1-vw}{1-v},\dfrac{1-vw}{1-w}).\]
\end{cor}

\subsection{Properties}
\label{sec:props}
From the definition of $Q'$, it is not difficult to describe its behaviour
under the polymatroid generalisation of many standard matroid operations;
we see that it retains versions of formulae true of the Tutte polynomial in many cases.
For instance, there is a polymatroid analogue of the direct sum of matroids:
given two polymatroids $M_1=(E_1,r_1),M_2=(E_2,r_2)$ with disjoint ground sets,
their \emph{direct sum} $M=(E,r)$ has ground set $E = E_1\sqcup E_2$
and rank function $r(S) = r_1(S\cap E_1)+r_2(S\cap E_2)$.

\begin{prop}Let $M_1\oplus M_2$ be the direct sum of two polymatroids $M_1$ and $M_2$.
Then 
\[Q'(M_1\oplus M_2;x,y)=\dfrac{Q'(M_1;x,y)Q'(M_2;x,y)}{x+y-1}.\]
\end{prop}

In particular, in the matroid setting where one of the summands is a loop or a coloop, we obtain:

\begin{cor}
Take a matroid $M=(E,r).$ Let $M'=M\cup\{e\}$ where $e$ is either a loop or a coloop. Then $Q'(M';x,y)=(x+y-1)Q'(M;x,y).$
\end{cor}

Next, $Q'$ exchanges its two variables under duality, as does the Tutte polynomial.
The best analogue of duality for polymatroids requires a parameter $s$
greater than or equal to the rank of any singleton; then if $M=(E,r)$ is a polymatroid,
its \emph{$s$-dual} is the polymatroid $M^*=(E,r^*)$ with 
\[r^*(S) = r(E) + s|E\setminus S| - r(E\setminus S).\]
The 1-dual of a matroid is its usual dual.

\begin{prop}
For any polymatroid $M=(E,r)$ and any $s$-dual $M^*$ of~$M$, 
$Q'(M^*;x,y)=Q'(M;y,x).$
\end{prop}

The invariant $Q'$ is a polytope valuation of polymatroids.  That is:
\begin{prop}
Let $\mathcal F$ be a polyhedral complex whose total space is a polymatroid base polytope
$P(M)$, and each of whose faces $F$ is a polymatroid base polytope $P(M(F))$.  Then
\[Q'(M;x,y) = \sum_{\mbox{\scriptsize $F$ a face of $\mathcal F$}}
(-1)^{\dim(P(M))-\dim F} Q'(M(F);x,y).\]
\end{prop}

For example, if $M$ is a matroid and we relax a circuit-hyperplane, we get the following result:

\begin{cor}
Take a matroid $M=(E,r)$ and let $C\subset E$ be a circuit-hyperplane of $M$. Let $M'$ be the matroid formed by relaxing $C$. Then $Q'(M;x,y)=Q'(M';x,y)-x^{n-r-1}y^{r-1}$.
\end{cor}

\section{Coefficients}
\label{sec:coeffs}
In this section, $M$ will always refer to a matroid, unless stated otherwise. Some coefficients of the Tutte polynomial provide structural information about the matroid in question. Let $b_{i,j}$ be the coefficient of $x^iy^j$ in $T(M;x,y)$. It is well known that $M$ is connected only if $b_{1,0}$, known as the \emph{beta invariant}, is non-zero; moreover, $b_{1,0} = b_{0,1}$. 
Not every coefficient yields such an appealing result, though of course
they do count the bases with internal and external activity of fixed sizes. 
We are able to provide a geometric interpretation of the coefficients of $Q'(M;x,y)$ when $M$ is a matroid, which is the focus of this section.

In order to do this, we will make use of a particular regular mixed subdivision of $u\Delta +P(M)+t\nabla$. 
This will be 
the regular subdivision determined by the ``lifted'' polytope $(P(M)\times \{0\})+\text{Conv}\{(u{\bf e}_i,\alpha_i)\}+\text{Conv}\{(-t{\bf e}_i,\beta_i)\}$
lying in $\mathbb R^E\times\mathbb R$, 
where $\alpha_1<\cdots <\alpha_n$, $\beta_1<\cdots <\beta_n$ are positive reals. 
When $t=u=1$, the associated height function
on the lattice points of $\Delta +P(M)+\nabla$
is 
\[h(x):=\text{min}\{\alpha_i+\beta_j \ | \ x-{\bf e}_i+{\bf e}_j\in \mathcal B_M\};\]
in general, one subtracts $t$ standard basis vectors and adds $u$ of them.

Let $\mathcal{F}$ be the set of lower faces of the lifted polytope which are visible from below, that is, facets maximising some linear function $a\cdot x$ where $a\in(\mathbb R^E\times\mathbb R)^*$ is a linear functional with last coordinate $a_{n+1}=-1$. 
For each such face $F\in\mathcal{F}$, let $\pi(F)$ be its projection back to $\mathbb{R}^n$. Now $\{\pi(F) \ | \ F\in\mathcal{F}\}$ is a regular subdivision of $u\Delta+P(M)+t\nabla$.

The structure of the face poset of this polyhedral subdivision 
does not depend on $t$ and $u$ as long as these are positive.
%As such, we will sometimes let $t=u=1$ where we only care about this coarse structure.

\begin{dfn}
A cell $F+G+H$ of the mixed subdivision $\mathcal F$ is a \emph{top degree face}
when $G$ is a vertex of $P(M)$ and there exists no cell $F+G'+H$ of~$\mathcal F$ where $G'\nsupseteq G$. 
\end{dfn}

%We will write top-degree faces as $u\Delta_{1\cup X}+{\bf e}_B+t\nabla_{1\cup Y}$,
%where in the service of readability we write $1$ instead of $\{1\}$. 
We can now state the main result of this section:

\begin{thm}
\label{coeffs}
Take the regular mixed subdivision of $u\Delta+P(M)+t\nabla$ as described above. We have that $|[x^iy^j]Q'|$ counts the cells $F+G+H$ of the mixed subdivision where $G$ is a vertex of $P(M)$ and there exists no cell $F+G'+H$ where $G'\nsupseteq G$ and $i=dim(F),j=dim(H)$. 
\end{thm}

The key fact in the proof is the following.

\begin{prop}\label{prop:top degree}
In the subdivision $\mathcal{F}$, 
each of the lattice points of $u\Delta+P(M)+t\nabla$ lies in a top degree face.
\end{prop}

To expose the combinatorial content of this proposition, we need to describe
the top degree faces more carefully.  All top degree faces are of dimension $|E|-1$,
and therefore have the form $u\Delta_{X}+{\bf e}_B+t\nabla_{Y}$,
where $X$ and $Y$ are subsets of~$E$ so that $\Delta_X$ and $\nabla_Y$
are transverse and of dimensions summing to $|E|-1$, 
which imply that $X\cup Y=E$ and $|X\cap Y|=1$.  In fact 
the conditions on the $\alpha$ and $\beta$ imply that $X\cap Y=\{1\}$.

\begin{lemma}\label{lem:B}
Take subsets $X$ and $Y$ of~$E$ with $X\cup Y=E$ and $X\cap Y=\{1\}$.
There is a unique basis $B$ such that 
$u\Delta_{1\cup X}+{\bf e}_B+t\nabla_{1\cup Y}$ is a top-degree face. 
It is the unique basis $B$ such that
no elements of $X$ are externally inactive and no elements of $Y$ are internally inactive with respect to $B$, with reversed order on $E$.
\end{lemma}

In particular, there are exactly $2^{|E|-1}$ top degree faces, one for each
valid choice of $X$ and~$Y$.
The basis $B$ can be found using the simplex algorithm for linear programming
on $P(M)$, applied to a linear functional constructed from the $\alpha$ and $\beta$
encoding the activity conditions.
This procedure can be completely combinatorialised, giving a way to start from a randomly chosen initial basis and make a sequence of exchanges which yields a unique output $B$ regardless of the input choice.

Now Proposition~\ref{prop:top degree} reads:

\begin{sublemma} Any $x\in (u\Delta+P(M)+t\nabla)\cap\mathbb{Z}^n$ on $\mathcal{F}$ is of the form $({\bf e}_{j_1},\alpha_1)+\cdots +({\bf e}_{j_u},\alpha_u)+({\bf e}_B,0)+(-{\bf e}_{i_1},\beta_1)+\cdots +(-{\bf e}_{i_t},\beta_t)$ such that there exists a partition $X\sqcup Y=[n]\backslash 1$ where every $i$ is in $1\cup Y$, every $j$ is in $1\cup X$, and the algorithm of Lemma~\ref{lem:B}, given $X,Y$ yields $B$.
\end{sublemma}

We then form a poset $P$ where the elements are the top degree faces and all nonempty intersections of sets of these, ordered by containment. 
This poset is a subposet of the face lattice of the $(|E|-1)$-dimensional cube
whose vertices correspond to the top degree faces. 
Proposition~\ref{prop:top degree} shows that every lattice point of 
$u\Delta+P(M)+t\nabla$ lies in at least one face in~$P$.
The total number of lattice points is given by inclusion-exclusion on
the function on~$P$ assigning to each element of~$P$ the number of lattice points
in that face. Carrying this proof through produces the change of variables
transforming $Q$ to $Q'$,
yielding the result of the theorem.

The above proof immediately yields the following result:

\begin{cor}
The signs of the coefficients of $Q'(M;x,y)$ are alternating.
\end{cor}

This is parallel (if opposite) to the Tutte polynomial, where the coefficients are all positive.  The coefficients of $Q'$, up to sign, have the combinatorial interpretation
of counting elements of~$P$ of form $u\Delta_X + {\bf e}_B + t\nabla_Y$
by the cardinalities of $X\setminus\{1\}$ and~$Y\setminus\{1\}$.
In particular the top degree faces are counted by the collection
of coefficients of~$Q'$ of top degree (hence the name),
and the degree $|E|-1$ terms of $Q'$ are always $(x+y)^{|E|-1}$.

The following two combinatorial lemmas are used to establish the necessary structure of~$P$.

\begin{lemma}
\label{meet}
Take two distinct partitions $(X_1,Y_1)$,$(X_2,Y_2)$ of $[n]\setminus\{1\}$. The algorithm finds two bases $B_1$, $B_2$ such that we have two top degree cells $T_i=\Delta_{1\cup X_i}+{\bf e}_{B_i}+\nabla_{1\cup Y_i}$, $i\in\{1,2\}$. If $T_1\cap T_2\neq\emptyset$, then $B_1=B_2$.
\end{lemma}

\begin{lemma}
Take two distinct partitions $P_1=(X_1,Y_1)$, $P_2=(X_2,Y_2)$ of $[n]\setminus\{1\}$ such that their corresponding top degree cells contain a common point $p$. Now let $P_3=(X_3,Y_3)$ be a partition of $[n]\setminus\{1\}$ such that if $x\in X_1\cap X_2$ then $x\in X_3$, and if $y\in Y_1\cap Y_2$ then $y\in Y_3$. Then $p\in T_3$, and
the bases found by Lemma~\ref{lem:B} are $B_3=B_1=B_2$.
\end{lemma}

The appearance of basis activities in Lemma~\ref{lem:B} reveals that
$P$ is intimately related to a familiar object in matroid theory,
the \emph{Dawson partition} \cite{dawson}. 
Give the lexicographic order to the power set $\mathcal P(E)$.
A partition of $\mathcal P(E)$ into intervals $[S_1,T_1],\ldots,[S_p,T_p]$ with indices such that $S_1<\ldots<S_p$  is a Dawson partition if and only if $T_1<\ldots < T_p$. 
Every matroid gives rise to a Dawson partition in which these intervals are $[B\setminus\mathrm{Int}(B),B\cup\mathrm{Ext}(B)]$ for all $B\in\mathcal B_M$.

\begin{prop}
Let $[S_1,T_1],\ldots,[S_p,T_p]$ be the Dawson partition of~$M$.
The poset $P$ is a disjoint union of face posets of cubes $C_1,\ldots,C_p$
where the vertices of $C_i$ are the top-degree faces $u\Delta_X+{\bf e}_B+t\nabla_Y$ 
such that $X\in[S_i,T_i]$.
\end{prop}

%We now look at the lifted polytope and want to pick out the top-degree faces --- take any point in such a face, and we need to be able to not ``move down'' to a point of lower height. 
%Let $X$ be the set of positions of $j$'s, and $Y$ be the position of $i's$, such that $X\sqcup Y=[n]\setminus \{1\}$.
%This requires all $j$'s in the above lemma to not be externally active, and all $i$'s to not be internally active, with respect to the basis $B$, and reversed ordering on $E$. 

Here is an example to illustrate this construction of $Q'$ 
and show that Theorem~\ref{coeffs} fails for polymatroids.
\begin{ex}
The left of Figure~\ref{fig:2} displays the subdivision $\mathcal F$ for the sum of Example~\ref{ex:1}.
\begin{figure}
\begin{center}
\includegraphics[scale=1]{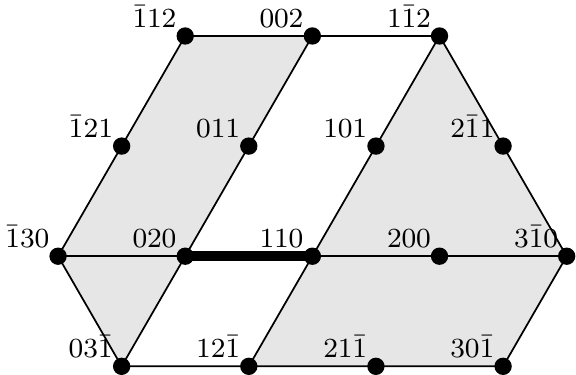}
\hspace{0.5in}
\includegraphics[scale=1]{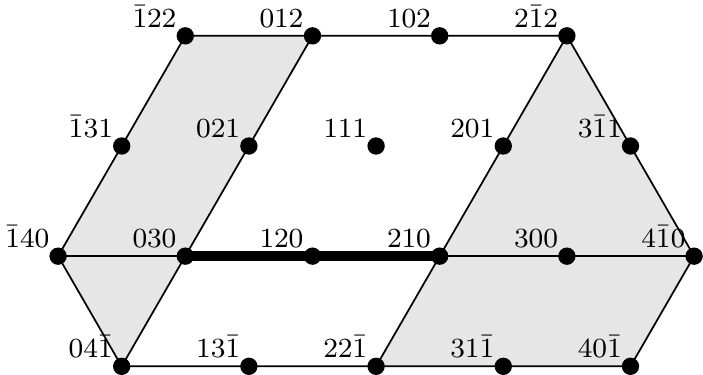}
\end{center}
\caption{At left, the regular subdivision $\mathcal F$ associated
to the Minkowski sum of Example~\ref{ex:1}, with 
$P(M)$ bolded and the top degree faces shaded in grey.
At right, the regular subdivision for a related polymatroid,
still with $(t,u)=(2,1)$.}\label{fig:2}
\end{figure}
We see that the four grey top degree faces contain all the lattice points between them,
and the poset $P$ contains two other faces which are pairwise intersections thereof,
the horizontal segment on the left with $(X,Y) = (1,12)$ and the one on the right with $(X,Y) = (12,1)$.
These are indeed enumerated, up to the alternation of sign,
by the polynomial $Q'(M;x,y) = x^2 + 2xy + y^2 - x - y$ found earlier.

By contrast, the right of the figure displays $\mathcal F$
for the polymatroid $M_2$ obtained by doubling the rank function of~$M$.
The corresponding polynomial is $Q'(M_2;x,y) = x^2 + 2xy + y^2 - 1$,
in which the signs are not alternating, dashing hopes of a similar enumerative interpretation.
In the figure we see that there are lattice points not on any grey face.
\end{ex}

\nocite{*}
\bibliographystyle{abbrvnat}
% use the following instead if you encounter problems 
%\bibliographystyle{alpha}
\bibliography{library}
\label{sec:biblio}

\end{document}